\documentclass{amsart}
\usepackage[utf8]{inputenc}
\usepackage{amsmath,amstext,amsopn,amssymb, amsfonts, amsthm, mathtools}
\usepackage{mathrsfs}
\usepackage[dvipsnames]{xcolor}
\usepackage[hmarginratio=1:1]{geometry}
\usepackage{bbm}
\usepackage{pdfsync}
\usepackage{tikz}
\usepackage{enumitem}

\usepackage[pdftex,bookmarks,colorlinks,breaklinks]{hyperref}  
\definecolor{dullmagenta}{rgb}{0.4,0,0.4}   
\definecolor{darkblue}{rgb}{0,0,0.4}
\definecolor{darkgreen}{rgb}{0,0.4,0}
\hypersetup{linkcolor=darkblue,citecolor=blue,filecolor=dullmagenta,urlcolor=darkblue} 

{}

\newtheorem*{definition*}{Definition}
\newtheorem{theorem}{Theorem}
\newtheorem*{theorem*}{Theorem}
\newtheorem{lemma}[theorem]{Lemma}
\newtheorem*{lemma*}{Lemma}

\newcommand{\dem}{\noindent {\bf Proof. }}
\newcommand{\demA}{\noindent {\bf Proof of Theorem \ref{thmDoublingCZ}. }}
\newcommand{\demB}{\noindent {\bf Proof of Theorem \ref{thmNonDoublingCZ}. }}

\newcommand{\fin}{\hspace*{\fill} $\square$ \vskip0.2cm}

\newtheorem{remark}[theorem]{Remark}
\newtheorem*{remark*}{Remark}

\newtheorem*{problem*}{Problem}
\numberwithin{equation}{section}
\numberwithin{theorem}{section}

\theoremstyle{theorem}
\newtheorem{ltheorem}{Theorem}

\newcommand{\G}{\mathrm{G}}

\newcommand{\Z}{\mathbb{Z}}
\newcommand{\R}{\mathbb{R}}
\newcommand{\Mb}{\mathbb{M}}

\newcommand{\Nb}{\mathbb{N}}

\newcommand{\A}{\mathcal{A}}
\newcommand{\V}{\mathcal{L}(\mathrm{G})}
\newcommand{\B}{\mathcal{B}}

\newcommand{\C}{\mathcal{C}}
\newcommand{\M}{\mathcal{M}}

\newcommand{\E}{\mathsf{E}}

\newcommand{\bd}{b_{\mathrm{d}}}
\newcommand{\boff}{b_{\mathrm{off}}}

\newcommand{\bdj}{b_{\mathrm{d},j}}
\newcommand{\boffj}{b_{\mathrm{off},j}}
\newcommand{\boffQ}{b_{\mathrm{off},Q}}

\newcommand{\1}{\mathbf{1}}

\newcommand\norm[1]{\left\| #1 \right\|}
\newcommand\md[1]{\left| #1 \right|}

\newcommand\set[1]{\left\lbrace #1 \right\rbrace}
\newcommand{\Ind}{1}
\newcommand{\supp}{\mathrm{supp}}
\makeatletter
\newcommand{\customlabel}[2]{%
   \protected@write \@auxout {}{\string \newlabel {#1}{{#2}{\thepage}{#2}{#1}{}} }%
   \hypertarget{#1}{#2}
}
\makeatother

\def\XXint#1#2#3{{\setbox0=\hbox{$#1{#2#3}{\int}$}
     \vcenter{\hbox{$#2#3$}}\kern-.5\wd0}}

\makeindex
\begin{document}

\addtolength{\parskip}{+1ex}

\null

\vskip-20pt 

\null

\title[Spectral multipliers in group algebras]{Spectral multipliers in group algebras \\ and noncommutative Calder\'on-Zygmund theory}

\author[Cadilhac, Conde-Alonso, Parcet]{L\'eonard Cadilhac, Jos\'e M. Conde-Alonso and Javier Parcet}

\begin{abstract}
In this paper we solve three problems in noncommutative harmonic analysis which are related to endpoint inequalities for singular integrals. In first place, we prove that an $L_2$-form of H\"ormander's kernel condition suffices for the weak type (1,1) of Calder\'on-Zygmund operators acting on matrix-valued functions. To that end, we introduce an improved CZ decomposition for martingale filtrations in von Neumann algebras, and apply a very simple unconventional argument which notably avoids pseudolocalization. In second place, we establish as well the weak $L_1$ endpoint for matrix-valued CZ operators over nondoubling measures of polynomial growth, in the line of the work of Tolsa and Nazarov/Treil/Volberg. The above results are valid for other von Neumann algebras and solve in the positive two open problems formulated in 2009. An even more interesting problem is the lack of $L_1$ endpoint inequalities for singular Fourier and Schur multipliers over nonabelian groups. Given a locally compact group $\G$ equipped with a conditionally negative length $\psi: \G \to \R_+$, we prove that Herz-Schur multipliers with symbol $m \circ \psi$ satisfying a Mikhlin condition in terms of the $\psi$-cocycle dimension are of weak type $(1,1)$. Our result extends to Fourier multipliers for amenable groups and impose sharp regularity conditions on the symbol. The proof crucially combines our new CZ methods with novel forms of recent transference techniques. This $L_1$ endpoint gives a very much expected inequality which complements the $L_\infty \to \mathrm{BMO}$ estimates proved in 2014 by Junge, Mei and Parcet.  
\end{abstract}

\allowdisplaybreaks

\maketitle


\vskip-20pt

\null

\section*{\bf Introduction} \label{section:introduction}

Fourier multipliers in group von Neumann algebras were early recognized by Haagerup as a key tool to reveal the geometry of these algebras via approximation properties \cite{Haa2,CanHaa,Haa1, HaaLaa}. The $L_p$ theory has provided in the last decade new rigidity theorems via Fourier and Schur approximation \cite{LaatSalle, LafforgueSalle, parcet-ricard-salle} and extensions of the H\"ormander-Mikhlin multiplier theorem with a broad interpretation of tangent spaces for nonabelian groups in terms of their cohomology \cite{glez-junge-parcetASENS, junge-mei-parcetSmooth, junge-mei-parcetJEMS, parcet-ricard-salle}. In 2010, at an early stage in the typing process of \cite{junge-mei-parcetSmooth}, Marius Junge already formulated the problem of establishing weak type $L_1$ bounds for Mikhlin multipliers in group von Neumann algebras. This natural question was appropriate since a noncommutative form of Calder\'on-Zygmund theory was launched in 2009, appearing as the next natural step after the development of operator space theory and noncommutative martingale inequalities. Unfortunately, the first contribution in this direction \cite{parcet2009} has been so far the only one to establish significantly new weak $L_1$ endpoint inequalities for noncommutative singular integrals. The main result in \cite{junge-mei-parcetSmooth} was the $L_\infty \to \mathrm{BMO}$ endpoint inequality for Mikhlin multipliers in discrete group von Neumann algebras. Both the Mikhlin condition and the BMO space are encoded in terms of a conditionally negative length $\psi: \G \to \R_+$. Before that, the $L_p$-boundedness of Fourier multipliers in group algebras was widely unexplored up to isolated contributions. Later connections with rigidity, geometric group theory or noncommutative geometry were found in \cite{glez-junge-parcetASENS,junge-mei-parcetJEMS,parcet-ricard-salle} with a crucial contribution of noncommutative CZ theory. In this paper,  we shall develop stronger Calder\'on-Zygmund methods to establish weak $L_1$ endpoints for Mikhlin multipliers and solve two further open problems in the process.

Beyond the strong connection with Fourier analysis in group von Neumann algebras, the theory of noncommutative CZ operators immediately showed a great versatility to solve problems in other areas. The main result in \cite{parcet2009} |CZ operators are of weak type $(1,1)$ when acting on Euclidean matrix-valued functions| was the key tool for the recent solution of the Nazarov-Peller conjecture \cite{CaspersetalAJM}, a strengthening in turn of the celebrated solution of the Krein conjecture \cite{PotapovSukochevActaMath} on operator Lipschitz estimates. In a different direction, pseudodifferential operators in noncommutative tori or the Heisenberg-Weyl algebra play a key role in Connes' noncommutative geometry. $L_p$-bounds for their commutative analogs are crucial in PDE and mathematical physics, but their proof requires CZ methods. This program was completed in the quantum realm for general quantum Euclidean spaces in \cite{glez-junge-parcetMAMS}, which involved the first CZ results for a fully noncommutative algebra. Other algebraic tools led to a CZ theory for general von Neumann algebras in \cite{junge-mei-parcet-xiaAIM}. The theory has developed rapidly in further directions over the last decade. Other interesting results and applications can be found in \cite{JungeSukochevetal,ChenXuZhiCMP, HongMeiT1, hong-lopezsanchez-martell-parcet2012, JungeRezvaniZengCMP, MeiParcetIMRN, XiaXiongXuAIM,XuTriebelMAMS} and the references therein. 

The proof of the main inequality in \cite{parcet2009} is extremely technical. After \cite{CaspersetalAJM}, it became a challenge to find a simpler argument which culminated in \cite{cadilhac2018}, where the original approach was significantly streamlined. In spite of that, the so-called pseudolocalization theorem was still unavoidable. This has been the fundamental obstruction to solve problems i) and ii) below, originally formulated in \cite[Remarks 5.5 and 5.6]{parcet2009}. In this paper, before considering Fourier and Schur multipliers, we shall need a new and much simpler approach to solve both of them:

\begin{itemize}
\item[\textbf{i)}] \textbf{Kernel regularity.} The inequality $$\sup_{y_1, y_2 \in \R^n} \int_{|x-y_1| > 2 |y_1 - y_2|} \big| k(x,y_1) - k(x,y_2) \big| \, dx < \infty$$ is known as H\"ormander's integral condition for the kernel $k(x,y)$. When the associated CZ operator $T_k$ is $L_2$-bounded, this condition is best known to imply the weak $L_1$ boundedness of $T_k$ and it can be found in almost any book on Euclidean harmonic analysis over the last 50 years. Compared to other regularity assumptions like the gradient or Lipschitz conditions it incorporates fundamental singular integral operators, prominently H\"ormander-Mikhlin multipliers. In the context of operator-valued functions, involved techniques from \cite{parcet2009} or even in the simpler proof from \cite{cadilhac2018} forced one to impose Lipschitz regularity on the CZ kernel. It was left as an open problem to decide whether it is possible to weaken the kernel regularity. We shall get an $L_2$-form of H\"ormander's condition in Theorem \ref{thmDoublingCZ} and we shall discuss the possible failure of Hormander's optimal condition in this setting. 

\vskip5pt

\item[\textbf{ii)}] \textbf{Nondoubling measures.} Calder\'on-Zygmund methods exploit the relation between metric and measure in the underlying space. This is particularly well understood when the measure $\mu$ is doubling. In other words, when there exists $\alpha>1$ nd $\beta >0$ such that $\mu(\alpha B) \le \beta \mu(B)$ for every ball $B$ in the given metric. Strong applications in geometric measure theory \cite{MasTolsa, NTVTb, TolsaPainleve, TolsaBilipschitz, TolsaGAFA} are intimately connected with a CZ theory for nondoubling measures of polynomial growth: $\mu(B)$ is just dominated up to a constant by a fixed power of its radius. It is definitely interesting to investigate a matrix-valued CZ theory for nondoubling measures  and explore its applications. All our attempts to prove nondoubling forms of pseudolocalization |the technique from \cite{cadilhac2018, parcet2009}| failed for several years and it seems highly difficult to decide whether it holds or not. We shall overcome this difficulty in Theorem \ref{thmNonDoublingCZ} with our new CZ tools and a nonregular filtration from \cite{condealonso-parcet2019}. 
\end{itemize}

Let $(\Omega, \mu)$ be either the $d$-dimensional Euclidean space with the Lebesgue measure or with other measure $\mu$ of polynomial growth, as described above. Let $T_k$ be a CZ operator with kernel $k$ defined on $(\Omega, \mu)$. We may extend it to matrix-valued functions $f: \Omega \to M_m$ by letting $T_k(f)$ be the matrix $(T_k(f_{ij}))$ where $f_{ij}(\omega)$ is the $(i,j)$-th entry of the $m \times m$ matrix $f(\omega)$. The first goal of this paper is to prove the inequality 
\begin{equation} \tag{W$L_1$} \label{GoalCZ}
\sup_{\lambda > 0} \, \lambda \int_\Omega \mathrm{tr} \big\{ |T_kf(\omega)| > \lambda \big\} \, d\mu(\omega) \le C_\Omega \int_\Omega \mathrm{tr} \big( |f(\omega)| \big) \, d\mu(\omega)
\end{equation}
for certain constant $C_\Omega$ independent of $m$ and kernels $k$ satisfying the weakest possible form of regularity. The strategy for \eqref{GoalCZ} in \cite{cadilhac2018, parcet2009} relied on a noncommutative martingale extension of CZ decomposition for scalar-valued functions. It looks as follows $$f \ = \ \overbrace{\sum_{j,k \in \widehat{\Z}} p_j f_{j \vee k} p_k}^{g_\lambda} + \overbrace{\sum_{j,k \in \widehat{\Z}} p_j (f - f_{j \vee k}) p_k}^{b_\lambda}.$$ Here, $\widehat{\Z} = \Z \cup \{\infty\}$ and $f_j = \mathsf{E}_j(f)$ is the standard $j$-th dyadic conditional expectation. Moreover, the $p_j$'s are a noncommutative analog of the sets where $f_j > \lambda$ but $f_{j-\ell} \le \lambda$ for every $\ell \in \Z_+$. In the commutative case $m=1$, if we consider the $\lambda$-level set for the dyadic maximal function $L_\lambda = \{M_df > \lambda\}$, $p_j$ is nothing but the union of maximal dyadic subcubes of $L_\lambda$ with side-length $2^{-j}$ and $p_\infty$ plays the role of $\R^d \setminus L_\lambda =  \{M_df \le \lambda \}$. 

Note that $p_j p_k = 0$ for $j \neq k$, but it does not kill the off-diagonal elements in the noncommutative case $m > 1$. The off-diagonal terms are harder to deal with than the diagonal ones. In the above decomposition, the off-diagonal terms for the good part $g_\lambda$ behave surprisingly worse than those for $b_\lambda$. For reasons we shall omit, the estimate of $T_k(g_{\lambda, \mathrm{off}})$ required to have a precise knowledge of the $L_2$-rate of decay of $T_kh$ away from $\supp \, h$. Roughly speaking, given a norm 1 function $h \in L_2(\R^d)$ with $\supp \, h = \Sigma$ compact, we need $$\int_{\R_+} \Big( \int_{\R^d \setminus \delta \Sigma} |T_k h(\omega)|^2 \, d\mu(\omega) \Big)^{\frac12} \, d\delta \, < \, \infty.$$ This pseudolocalization principle shows that the $L_2$-mass of $T_kh$ is somehow concentrated around $\supp \, h$ and is a fundamental obstruction to solve the above mentioned problems. Pseudolocalization was succesfully proved in \cite{parcet2009} with exponential decay |\cite{cadilhac2018} gives a simpler proof| for the Lebesgue measure (or any other doubling mesure) under the Lipschitz kernel condition, much stronger than the H\"ormander integral condition. As explained in \cite{parcet2009}, the relation between pseudolocalization and the $T1$ theorem |for which one needs Lipschitz regularity or similar| makes impossible to weaken the kernel assumptions substantially when following that approach. In addition, it is still unknown whether pseudolocalization holds for nondoubling measures. 

Our first contribution is a new CZ decomposition which circumvents the pseudolocalization technique and still yields $L_1$ endpoints for matrix-valued functions. The solution for doubling measures is surprisingly simple, just replace $j \vee k = \max \{j,k\}$ by the minimum $j \wedge k$. The key advantage is that $p_j f_{j \wedge k} p_k = 0$ for $j \neq k$: 
\begin{equation} \tag{CZD} \label{CZdec}
f \ = \ \underbrace{qfq + \sum_{j \in \Z} p_j f_j p_j}_{\mathrm{New} \ g_\lambda} + \underbrace{\frac12 \sum_{j \in \Z} p_j (f - f_j) (q_{j-1} + q_j) + (q_{j-1} + q_j) (f - f_j) p_j}_{\mathrm{New} \ b_\lambda},
\end{equation}
where $q_j$ is a noncommutative form of the set of points whose dyadic ancestor $Q$ of side-length $2^{-j}$ satisfies $Q \nsubseteq L_\lambda$. This eliminates off-diagonal terms in $g_\lambda$ and thereby pseudolocalization. The difficulty |faced in the preparation of \cite{parcet2009}, but never solved| is to prove \eqref{GoalCZ} using this new form of the bad part $b_\lambda$. It does not follow from canonical adaptations of classical arguments and requires an unconventional new idea which works for an $L_2$-form of Hormander's condition.   

\begin{ltheorem} \label{thmDoublingCZ}
Let $T_k$ be an $L_2$-bounded \emph{CZ} operator on $(\R^d,dx)$ with $$\sup_{Q \ \mathrm{dyadic}} \sum_{j \ge 1} \, \sup_{y \in Q} \Big( 2^{jd} \ell(Q)^d \int_{2^j \ell(Q) \le |x-c_Q| \le 2^{j+1} \ell(Q)} \big| k(x,y) - k(x,c_Q) \big|^2 \, dx \Big)^{\frac12} < \infty.$$ Then, the inequality \emph{\eqref{GoalCZ}} holds up to a dimensional constant $C_d$ which is independent of $m$. 
\end{ltheorem}

Our kernel condition is considerably weaker than Lipschitz regularity and close enough to H\"ormander's condition so as to include the class of H\"ormander-Mikhlin multipliers, important singular integrals excluded in \cite{cadilhac2018, parcet2009}. Our proof is also much simpler. A slightly less flexible kernel condition has been independently found and announced in \cite{hong-lai-xu2020},  but the proof was omitted since their argument is parallel to ours and based on our techniques: \eqref{CZdec} and the unconventional estimate alluded above, both of which we communicated to the authors. (They also used these ideas to obtain interesting results for maximal singular integrals.) As we shall discuss, it is very unclear whether one can replace our condition by H\"ormander's integral condition. 

The nonhomogeneous problem entails additional difficulties. First, the classical arguments for $L_1$ endpoints of nondoubling CZ operators rely either on approximation of measures \cite{NTV} or on a suitable \lq centered\rq${}$ CZ decomposition \cite{tolsaWeakL1}, as opposed to dyadic. There are relevant obstructions to transfer either of them to the noncommutative setting. Second, although the former can be remedied by using the dyadic-like construction given in \cite{condealonso-parcet2019},  pseudolocalization arguments are not available for nondoubling measures. And even if they were, the approach from \cite{cadilhac2018, parcet2009} requires to take centered dilations of cubes which respect the measurability in the dyadic filtration. This cannot work in the nondoubling setting because of the necessary lack of regularity of any filtration that one can construct. We solve this with a CZ decomposition for nonregular filtrations. 

Our decomposition in the nondoubling framework necessarily deviates from \eqref{CZdec} and will be made explicit in the body of the paper. We should recall though some terminology and results from \cite{condealonso-parcet2019}, to be amplified later on. Let $(\Omega,\mu)$ be $\R^d$ equipped with a measure $\mu$ of $n$-polynomial growth: $\mu(B(x,r)) \le C_\mu r^n$ for $\mu$-almost every $x \in \supp(\mu)$ and every ball $B(x,r)$ centered at $x$ with radius $r$. A ball $B$ is called $(\alpha,\beta)$-doubling when $\mu(\alpha B) \le \beta \mu(B)$. The main result in \cite{condealonso-parcet2019} is the construction of a weak-$*$ dense two-sided (nonregular) martingale filtration $(\Sigma_k)_{k \in \Z}$ of atomic $\sigma$-algebras of $\supp (\mu)$ enjoying several key properties. One of these is the existence of $\alpha, \beta \in \R_+$ such that, for every atom $Q$ in $\bigcup_k \Sigma_k$, there exists an $(\alpha,\beta)$-doubling ball $B_Q$ which is comparable to the atom: $B_Q \subset Q \subset \alpha B_Q$. Theorem \ref{thmNonDoublingCZ} below solves problem ii).  

\begin{ltheorem} \label{thmNonDoublingCZ}
Let $T_k$ be an $L_2$-bounded \emph{CZ} operator on $(\R^d,\mu)$ with 
\begin{itemize}
\item $\displaystyle |k(x,y)| \, \lesssim \, \frac{1}{|x-y|^n}$,

\item $\displaystyle \sup_{Q \ \mathrm{atom}} \, \sum_{j \ge 1} \, \sup_{y \in B_Q} \, \Big( \mu(\C_{Q,j}) \int_{\C_{Q,j}} \big| k(x,y)-k(x,c_{B_Q}) \big|^2 \, d\mu(x) \Big)^{\frac12} \, < \, \infty$. 
\end{itemize}
for $\C_{Q,j} = \{2^j r(B_Q) \le |x -  c_{B_Q}| \le  2^{j+1} r(B_Q)\}$. Then, \emph{\eqref{GoalCZ}} holds with $C_{\Omega}$ independent of $m$. 
\end{ltheorem}

We now switch to Fourier and Schur multipliers. Let $\G$ be a locally compact unimodular group with Haar measure $\mu$. Let $\lambda$ be the left-regular representation of $\G$. In other words, $\lambda(g)$ is the unitary $(\lambda(g) \varphi)(h) = \varphi(g^{-1}h)$ on $L_2(\G,\mu)$. The group von Neumann algebra $\V$ is the weak-$*$ closure in $\B(L_2(\G))$ of the linear $\mathrm{span} \, \{ \lambda(g): g \in \G \}$. It is also the weak-$*$ closure of elements of the form $$f = \int_\G \widehat{f}(g) \lambda(g) \, d\mu(g) \quad \mbox{with} \quad \widehat{f} \in \C_c(\G).$$ The canonical trace $\tau$ on $\V$ is densely defined by $\tau(f) = \widehat{f}(e)$ for smooth enough $\widehat{f}$. This allows one to define the corresponding noncommutative $L_p$ and weak $L_p$ spaces on $\V$. We refer to \cite{junge-mei-parcetJEMS, PisierXuNCLp} for precise definitions. If $\G$ is abelian, recall that $L_p(\V)$ is just the natural $L_p$ space on the dual group of $\G$. Given a bounded measurable symbol $M: \G \to \mathbb{C}$, the Fourier multiplier associated to it is the linear map which intertwines multiplication via the Fourier transform $$T_M f = \int_\G M(g) \widehat{f}(g) \, d\mu(g) \quad \mbox{so that} \quad \widehat{T_Mf}(g) = M(g) \widehat{f}(g).$$

According to Schoenberg's theorem, a length $\psi: \G \to \R_+$ is conditionally negative exactly when the map $\lambda(g) \mapsto \exp(- t \psi(g)) \lambda(g)$ defines a Markov semigroup on $\V$. Its infinitesimal generator may be understood as a Laplacian $A_\psi: \lambda(g) \mapsto \psi(g) \lambda(g)$ on the group algebra of $\G$. Spectral multipliers are Fourier multipliers associated to symbols of the form $m \circ \psi$ for a conditionally negative length $\psi$ and some $m: \R_+ \to \mathbb{C}$. This class of multipliers is specially relevant since they arise by functional calculus on the Laplacian $A_\psi$ $$T_{m \circ \psi} = m(A_\psi).$$ H\"ormander-Mikhlin criteria for spectral multipliers have been deeply investigated in the remarkable work of Cowling, M\"uller, Ricci, Stein and their coauthors on nilpotent groups \cite{MartiniMuller,MullerCrelle,MullerRicciStein,RicciAJM}. Here we take the dual approach initiated by Haagerup and investigate analogous criteria on group von Neumann algebras. The information carried by conditionally negative lengths is rich enough to link Fourier multipliers in group algebras with Euclidean harmonic analysis. The main discovery in \cite{junge-mei-parcetSmooth} was that a H\"ormander-Mikhlin theory in group algebras can be deduced from Euclidean harmonic analysis, noncommutative Calder\'on-Zygmund theory and basic group cohomology. This line has been further investigated in \cite{glez-junge-parcetASENS, junge-mei-parcetJEMS, parcet-ricard-salle}, but none of these works provides any information on weak type $L_1$ estimates for Mikhlin multipliers in group von Neumann algebras. On the other hand, let $S_p(\G)$ be the Schatten $p$-class of compact operators on the Hilbert space $L_2(\G,\mu)$. Given a Fourier symbol $M: \G \to \mathbb{C}$, the Herz-Schur multiplier $S_M$ is the linear map on $S_2(\G)$ $$S_M(A) = \Big( M(gh^{-1}) A_{gh} \Big).$$

\begin{ltheorem} \label{thmFourierMult}
Let $\G$ be a locally compact group equipped with a $n$-dimensional conditionally negative length $\psi: \G \to \R_+$. Let $m: \R_+ \to \mathbb{C}$ be a Fourier symbol satisfying the Mikhlin condition in dimension $n$
$$\Big|\frac{d^k}{d \xi^k} m(\xi) \Big| \le C_{\mathrm{hm}} |\xi|^{-k} \quad \mbox{for all} \quad \xi \neq 0 \quad \mbox{and} \quad 0 \le k \le \Big[ \frac{n}{2} \Big] + 1.$$
Then, the following weak type $(1,1)$ inequalities hold$\hskip1pt :$ 
\begin{itemize}
\item[\emph{i)}] $S_{m \circ \psi}$ extends to a bounded map $S_{m \circ \psi}: S_1(\G) \to S_{1,\infty}(\G)$. 

\item[\emph{ii)}] $T_{m \circ \psi}: L_1(\V) \to L_{1,\infty}(\V)$ is also bounded for $\G$ discrete and amenable.
\end{itemize}
In fact, $S_{m \circ \psi} \otimes \mathrm{Id}_{\Mb_k}$ and $T_{m \circ \psi} \otimes \mathrm{Id}_{\Mb_k}$ remain $L_1 \to L_{1,\infty}$ bounded after $(k \times k)$-matrix amplification.
\end{ltheorem}

The dimension of the conditionally negative length $\psi$ is the dimension of the only cocycle $\beta$ for which $\psi(g) = |\beta(g)|^2$. Precise definitions will be given in the body of the paper. Theorem \ref{thmFourierMult} establishes a very much expected $L_1$ endpoint inequality for spectral multipliers in group von Neumann algebras, with sharp regularity order in terms of $\dim \psi$ as it follows from $\G = \R^n$ and $\psi(\xi) = |\xi|^2$ \cite[Theorem C]{junge-mei-parcetSmooth}. The proof presented here requires weak $L_1$ forms of some results of independent interest, which include transference between Fourier and Schur multipliers as well as their relation with twisted Fourier multipliers. It also involves a generalization of de Leeuw's approximation method together with our new CZ type results described above. Theorem \ref{thmFourierMult} also yields an $L_1$ endpoint for the Littlewood-Paley theorem in matrix and group algebras.  

The rest of the paper is divided into three sections. In Section \ref{section:CZdecomposition}, we shall present the new CZ decompositions for regular and nonregular martingale filtrations.  The proofs of Theorems \ref{thmDoublingCZ} and \ref{thmNonDoublingCZ} will be presented in Section \ref{section:Weak11}, while Fourier and Schur multipliers will be investigated in Section \ref{sectionFourier}.
     

\section{\bf Calder\'on-Zygmund decompositions} \label{section:CZdecomposition}


A noncommutative measure space $(\M,\tau)$ is a pair formed by a semifinite von Neumann algebra $\M$ equipped with a normal faithful normal trace $\tau$. The reader unfamiliar with noncommutative integration may replace $\M$ by the matrix algebra $\B(\ell_2)$ of bounded linear operators in $\ell_2$ with its natural trace $\mathrm{tr}$. A martingale filtration $\{\M_j\}_{j \in \Nb}$ of $\M$ is a nested family of von Neumann subalgebras $\M_j$ whose union is weak-$*$ dense in $\M$. Assume that for every $j \geq 1$, there is a normal conditional expectation $\E_j: \M \to \M_j$. A noncommutative martingale with respect to the above filtration is a sequence $(f_j)_{j \ge 1}$ in $\M$ such that $\mathsf{E}_j(f_k) = f_j$ for all $j \le k$. We refer to \cite{PisierMartingales} for an overview of the theory of noncommutative $L_p$-martingale inequalities. We will just need the so-called Cuculescu construction from \cite{cuculescu1971}. It is a noncommutative form of the weak type $(1,1)$ inequality for the Doob maximal function. Namely, given $f\in L_1(\M)_+$ and $\lambda>0$ there exists a nonincreasing sequence of projections $\{q_j = q_j(f,\lambda) \}_{j \in\Nb}$ with $q_0 = \mathbf{1}_\M$ such that the following properties hold:
\begin{itemize}
\item[\bf A.] $q_j \in \M_j$ for each $j \in \Nb$. \vskip1pt
\item[\bf B.] $q_j\E_j(f)q_j \leq \lambda q_j$ for each $j\in\Nb$. \vskip1pt
\item[\bf C.] $q_j$ commutes with $q_{j-1}\E_j(f)q_{j-1}$ for each $j\in\Nb$.
\item[\bf D.] If $q = \displaystyle \bigwedge_{j\in\Nb} q_j$, then $qfq \leq \lambda q$ and $\lambda \hskip1pt \tau(\mathbf{1}_\M-q) \leq \|f\|_1$.
\end{itemize} 
The proof of Cuculescu's theorem was originally formulated over noncommutative probability spaces, but the same argument applies for non necessarily finite von Neumann algebras. We refer to \cite{parcet2009,ParRan} for a more in depth discussion. 


\subsection{Regular martingales} 

The martingale filtrations above are called regular when there exists an absolute constant $c_{\mathrm{reg}} > 0$ satisfying $\E_j f \leq c_{\mathrm{reg}} \, \E_{j-1} f$ for every $f \in \M_+$. We may now introduce our new decomposition of CZ type for a regular filtration $\{\M_j\}_{j \in \Nb}$. Namely, given $f \in L_1(\M)_+$ and $\lambda >0$, consider the Cucuclescu's projections $q_j = q_j(f,\lambda)$. Define $p_j = q_{j-1} - q_j \in \M_j$ for $j \ge 1$, so that 
$$\sum_{j \ge 1} p_j = \mathbf{1}_\M - q.$$ Writing $f_j$ for $\E_j(f)$, this readily gives 
\begin{eqnarray*}
f & = & qfq + \sum_{j,k \ge 1}  p_j f p_k + \sum_{j \ge 1} p_j f q + \sum_{k \ge 1} q f p_k \\ & = & qfq + \sum_{j,k \ge 1}  p_j f_{j \wedge k} p_k + \sum_{j \ge 1} p_j f_j q + \sum_{k \ge 1} q f_k p_k \\ & + & \sum_{j,k \ge 1}  p_j (f - f_{j \wedge k}) p_k + \sum_{j \ge 1} p_j (f - f_j) q + \sum_{k \ge 1} q (f - f_k) p_k.
\end{eqnarray*} 
Cuculescu's commutation relations give 
\begin{itemize}
\item $p_j f_j q = p_j q_{j-1} f_j q_{j-1} q = 0$ for all $j \ge 1$. 

\item $p_j f_{j \wedge k} p_k = p_j q_{j \wedge k - 1} f_{j \wedge k} q_{j \wedge k -1} p_k = 0$ for all $j \neq k$.
\end{itemize}
In particular, the above decomposition of $f$ simplifies as follows
\begin{eqnarray} \label{CZRegular} 
f & = & qfq + \sum_{j \ge 1} p_j f_j p_j + \sum_{j \ge 1} p_j (f - f_j) p_j \\ \nonumber & + & \sum_{j \ge 1} p_j (f - f_j) \Big[ q + \sum_{k > j} p_k \Big] + \sum_{k \ge 1} \Big[ q + \sum_{j > k} p_j \Big] (f - f_k) p_k \\ \nonumber & = & \underbrace{qfq + \sum_{j \ge 1} p_j f_j p_j}_{g} + \underbrace{\underbrace{\sum_{j \ge 1} p_j (f - f_j) p_j}_{b_\mathrm{d} = \sum_j b_{\mathrm{d},j}} + \underbrace{\sum_{j \ge 1} p_j (f - f_j) q_j + q_j (f - f_j) p_j}_{b_{\mathrm{off}} = \sum_j b_{\mathrm{off},j}}}_{b}.     
\end{eqnarray}
A straightforward reordering gives the decomposition \eqref{CZdec} alluded in the Introduction. The above decomposition recovers the standard dyadic CZ decomposition for classical functions, since $b_{\mathrm{off}}$ vanishes in that case by commutativity. Note as well that $b_{\mathrm{off}} = \sum_{j \ge 1} p_j f q_j + q_j f p_j$ by Cuculescu's commutation relations. 

\begin{lemma} \label{thm:CZdecRegular}
The regular decomposition \eqref{CZRegular} satisfies$\hskip1pt :$
\begin{itemize}
\item $\|g\|_1 \le \|f\|_1$ and $\|g\|_\infty  \leq c_{\mathrm{reg}} \lambda$, 

\vskip7pt

\item $\displaystyle \sum_{j \ge 1} \|\bdj\|_1 \leq 2 \|f\|_1$ and $\displaystyle \E_j(\bdj) = \E_j(\boffj) =0$.
\end{itemize}
\end{lemma} 

\dem Compared to \cite{cadilhac2018, parcet2009}, the only novelty here is the mean 0 of $b_{\mathrm{off},j}$, which is trivial.
\fin


\subsection{Nonregular martingales} We now turn our attention to nonregular martingales. The main difficulty here is the lack of an $L_\infty$ estimate for the good part of the function, which requires regularity. The solution, as first introduced in \cite{lopezsanchez-martell-parcet2014} and then adapted to the study of CZ operators in \cite{condealonso-parcet2019}, is to modify the regular CZ decomposition by adding and subtracting $\E_{j-1}(p_jfp_j)$ instead of $p_jf_jq_j$. The resulting decomposition is \vskip-10pt
\begin{equation} \label{CZNonRegular} 
f \ = \ \underbrace{qfq + \sum_{j \ge 1} \E_{j-1}(p_j f p_j)}_{g} + \underbrace{\underbrace{\sum_{j \ge 1} p_j f p_j - \E_{j-1}(p_j f p_j)}_{b_\mathrm{d} = \sum_j b_{\mathrm{d},j}} + \underbrace{\sum_{j \ge 1} p_j f q_j + q_j f p_j}_{b_{\mathrm{off}} = \sum_j b_{\mathrm{off},j}}}_{b}.     
\end{equation}
This does not yield a uniform estimate, but an $L_2$ substitute is enough for most purposes. 

\begin{lemma} \label{thm:CZdecNonRegular}
The nonregular decomposition \eqref{CZNonRegular} satisfies$\hskip1pt :$
\begin{itemize}
\item $\|g\|_1 \le \|f\|_1$ and $\|g\|_2^2  \leq 6 \lambda \|f\|_1$, 

\vskip7pt

\item $\displaystyle \sum_{j \ge 1} \|\bdj\|_1 \leq 2 \|f\|_1$ and $\displaystyle \E_j(\bdj) = \E_j(\boffj) =0$.
\end{itemize}
\end{lemma} 

\dem The results for $\bd$ and $\boff$ are identical to the ones in the regular case, because regularity plays no role in their proof. Since the $L_1$-estimate for $g$ is trivial, we are left with the $L_2$ estimate. Since $q_j$ commutes with $q_{j-1}f_jq_{j-1}$, 
$$\E_{j-1}(p_j f p_j) = \E_{j-1}(p_jf_jp_j) = q_{j-1}f_{j-1}q_{j-1} - \E_{j-1}(q_jf_jq_j).$$
In particular, we get
$$\Big\| \sum_{j \ge 1} \E_{j-1}(p_j f_j p_j) \Big\|_2^2 \le 2 \Bigg( \underbrace{\Big\| \sum_{j \ge 1} q_j f_j q_j - \E_{j-1}(q_j f_j q_j) \Big\|_2^2}_{\mathrm{A}} + \underbrace{\Big\|\sum_{j \ge 1} q_j f_j q_j - q_{j-1}f_{j-1}q_{j-1} \Big\|_2^2}_{\mathrm{B}} \Bigg).$$
As it was proved in \cite[Lemma 3.4]{randrianantoanina2002}, we have that
$$\big\| q_jf_jq_j - \mathsf{E}_{j-1}(q_jf_jq_j) \big\|_2^2 \le 2 \Big( \big\| q_jf_jq_j \big\|_2^2 - \big\|q_{j-1}f_{j-1}q_{j-1}\|_2^2 \Big) + 6 \lambda \tau (q_{j-1}f_{j-1}q_{j-1} - q_{j}f_{j}q_{j}).$$
Therefore, by orthogonality of martingale differences, summation over $j$ gives
$$\mathrm{A} = \sum_{j \ge 1} \big\| q_j f_jq_j - \E_{j-1}(q_{j}f_jq_{j}) \big\|_2^2 \le 2 \big\| qfq \big\|_2^2 - 6\lambda\tau(qfq) \le 2\lambda \|f\|_1$$ since $qfq \leq \lambda q$ from \cite[Section 4.1]{parcet2009}. The telescopic sum in $\mathrm{B}$ easily yields that $\mathrm{B} \le \lambda \|f\|_1$. \fin

\begin{remark}
\emph{In the next sections we will derive consequences of the decompositions introduced above that will require some work. There is however an immediate application of Lemma \ref{thm:CZdecNonRegular} which is worth mentioning. Operators adapted to the underlying martingale structure such as the dyadic Hilbert transform and other Haar shifts can be shown to be of weak type $(1,1)$ |when acting on matrix-valued functions, with constants independent of the matrix size| under the natural assumptions on their symbols, as considered in \cite{condealonso-lopezsanchez}. The new decomposition makes the proof of the weak type estimates almost trivial for these operators. We omit the details.}
\end{remark}


\section{\bf Noncommutative Calder\'on-Zygmund operators} \label{section:Weak11}

Let $\mu$ be a locally finite Radon measure on $\Omega = \R^d$ and consider a noncommutative measure space $(\M,\tau)$. Let us construct the tensor product von Neumann algebra $\A = L_\infty(\Omega,\mu) \bar\otimes \M$ with its natural trace $$\varphi(f) = \int_\Omega \tau(f(\omega)) \, d\mu(\omega).$$ Elements in $\A$ may be identified with functions $f: \Omega \to \M$ such that $\omega \mapsto \|f(\omega)\|_\M$ is essentially bounded. Similarly, the space $L_p(\A)$ becomes the space of $p$-integrable vector-valued functions $f: \Omega \to L_p(\M)$ with values in the noncommutative $L_p$ space associated to $\M$. As we did in the Introduction, we shall fix $\M = M_m$ with its canonical trace for concreteness, though all of our new results apply equally well to any noncommutative measure space $(\M,\tau)$. 

In this section, we shall consider kernels $k: \R^{2d} \setminus \{(x,y):x=y\} \to \mathbb{C}$ satisfying the regularity conditions imposed in Theorems \ref{thmDoublingCZ} or \ref{thmNonDoublingCZ}. Accordingly, we consider linear maps which admit the following integral representation for matrix-valued functions $f: \Omega \to M_m$ $$T_kf(x) = \int_{\Omega} k(x,y)f(y) \, d\mu(y) \quad \mbox{for every}\quad x \notin \supp_{\R^d} (f).$$ Here we write $\supp_{\R^d} (f) = \supp \hskip1pt \|f\|_{M_m}$ to denote the Euclidean support of the function $f$, as opposed to its support as an operator affiliated to $\A$. We shall also consider two-sided martingale filtrations $\{\A_j\}_{j \in \Z}$ of the form 
$$
\A_j = L_\infty(\Omega,\Sigma_j,\mu) \otimes (M_m,\mathrm{tr}),
$$
where $\{ \Sigma_j \}$ is a nested family of atomic $\sigma$-algebras on $\Omega = \R^d$. We shall write $\Pi(\Sigma_j)$ and $\Pi(\Sigma)$ to denote the set of atoms in $\Sigma_j$ and $\bigcup_j \Sigma_j$ respectively. Eventually, we shall apply CZ decompositions at a fixed height $\lambda$ to elements $f \in L_1(\A)_+$ with respect to the above filtration. By density considerations as in \cite[Section 3]{parcet2009}, we may assume that $\supp_{\R^d}(f)$ is compact and Cuculescu's projections $q_j = \mathbf{1}_\A$ for all $j \le m_{f,\lambda}$ and some $m_{f,\lambda} \in \Z$. In other words, for fixed $(f,\lambda)$ we will not use the full two-sided filtration, but  a truncation from $m_{f,\lambda}$. In fact, we may assume for simplicity that $m_{f,\lambda} = 0$. This allows us to use Cuculescu's construction in what follows as a black box. Due to the atomic nature of our filtrations, Cuculescu's projection $q_j=q_j(f,\lambda)$ must be constant on each atom $Q \in \Pi(\Sigma_j)$ and there must exists projections $\pi_Q, \xi_Q$ in $M_m$ such that  
\begin{eqnarray*}
q_j \!\!\! & = & \!\!\! \sum_{Q\in\Pi(\Sigma_j)} q_Q \ = \sum_{Q\in\Pi(\Sigma_j)} \Ind_Q \otimes \xi_Q, \\
p_j \!\!\! & =& \!\!\! \sum_{Q\in\Pi(\Sigma_j)} p_Q \ = \sum_{Q\in\Pi(\Sigma_j)} \Ind_Q \otimes \pi_Q.
\end{eqnarray*}



\subsection{Doubling measures} Let $\mu$ be the $d$-dimensional Lebesgue measure. As will be evident from the proof, all our estimates hold whenever $\mu$ is any doubling measure. Let $T_k$ be a CZ operator with kernel $k$ as above. The filtration $\{\A_j\}_{j \in\Z}$ that we use is the same that is used in the scalar valued case: the one generated by the dyadic system $\mathscr{D}$. This can be defined as
$$
\mathscr{D} = \bigcup_{j\in\Z} \mathscr{D}_j = \bigcup_{j\in\Z} \Big\{ 2^{-j} \cdot [k_1,k_1+1) \times [k_2,k_2+1) \times \ldots [k_d,k_d+1): (k_1, \ldots, k_d) \in \Z^d \Big\}.
$$
Let $\Sigma_j = \sigma(\mathscr{D}_j)$ so that $\Pi(\Sigma_j) = \mathscr{D}_j$. By the translation invariance of the Lebesgue measure, the dyadic filtration is regular with constant $2^d$. Given $f \in L_1(\A)_+$ compactly supported and $\lambda >0$, we consider Cuculescu's projections $q_j, p_j$ and $\xi_Q, \pi_Q$ for $Q \in \Pi(\Sigma_j)$. Define 
$$
\zeta := \mathbf{1}_\A - \bigvee_{j \ge 1} \bigvee_{Q \in \mathscr{D}_j} 1_{5Q} \pi_Q. 
$$
As explained in \cite{cadilhac2018, parcet2009}, the projection $\mathbf{1}_\A - \zeta$ represents an Euclidean dilation of the CZ maximal cubes of $f$ at height $\lambda$, as will become clear from its role in the proof of the weak boundedness of $T_k$. We will need the following properties, whose proof can be found in \cite[Lemma 3.4]{cadilhac2018}.

\begin{lemma}\label{lem:zetaproydoubling} The projection $\zeta$ satisfies
$$
\varphi(\mathbf{1}_\A - \zeta) \lesssim \frac{\|f\|_{L_1(\A)}}{\lambda}.
$$
Moreover, we have $\zeta(x)p_j(y) = p_j(y)\zeta(x) = 0$ whenever $y \in Q \in \mathscr{D}_j$ and $x \in 5Q$.
\end{lemma}


\demA It suffices to show that
\begin{equation}\label{eq:goal1}
\varphi \big\{ |T_kf|>7\lambda \big\} \, \lesssim \, \frac{\|f\|_{L_1(\A)}}{\lambda}.
\end{equation}
for $f \in L_1(\A)_+$ with $\supp_{\R^d}(f)$ compact. We write 
\begin{eqnarray*}
T_kf & = & T_kg \\ & + & (\mathbf{1}-\zeta)T_k\bd + \zeta T_k\bd(\mathbf{1}-\zeta) + \zeta T_k\bd \zeta \\ & + & (\mathbf{1}-\zeta)T_k\boff + \zeta T_k\boff(\mathbf{1}-\zeta) + \zeta T_k\boff \zeta.
\end{eqnarray*}
where $f = g + \bd + \boff$ using Lemma \ref{thm:CZdecRegular} at height $\lambda$. Therefore, we can estimate
\begin{eqnarray*}
\varphi(\set{|T_kf|>7 \lambda}) & \le & \varphi \big\{ |T_kg|>\lambda \big\} \\
& + & \varphi \big\{ |(\1-\zeta)T_k\bd|>\lambda \big\} + \varphi \big\{ |\zeta T_k\bd(\1-\zeta)|>\lambda \big\} + \varphi \big\{ |\zeta T_k\bd \zeta|>\lambda \big\} \\ & + & \varphi \big\{ |(\1-\zeta)T_k\boff|>\lambda \big\} + \varphi \big\{ |\zeta T_k\boff(\1-\zeta)|>\lambda \big\} + \varphi \big\{ |\zeta T_k\boff \zeta|>\lambda \big\} \\ & \lesssim & \varphi \big\{ |T_kg|>\lambda \big\} + \varphi \big\{ |\zeta T_k\bd \zeta|>\lambda \big\} + \varphi \big\{ |\zeta T_k\boff \zeta|>\lambda \big\} + \lambda^{-1} \|f\|_1
\end{eqnarray*}
by Lemma \ref{lem:zetaproydoubling} and Murray--von Neumann equivalence. The terms $g$ and $b_\mathrm{d}$ can be dealt with in a very similar way than in \cite{cadilhac2018, parcet2009}, just being careful enough to use H\"ormander's kernel condition  instead of Lipschitz regularity for $b_{\mathrm{d}}$. We include the details for the sake of completeness. 

\noindent \textbf{A. Standard estimates.} 
We clearly have 
$$
\varphi \big\{ |T_kg|>\lambda \big\} \leq \frac{1}{\lambda^2} \|T_kg\|_2^2 \lesssim \frac{1}{\lambda^2} \|g\|_2^2 \lesssim \frac{\|g\|_1}{\lambda} \leq \frac{\|f\|_1}{\lambda}.
$$
as a consequence of Chebychev inequality, the $L_2$ boundedness of $T_k$ and the estimate $\|g\|_\infty \leq 2^d \lambda$ from Lemma \ref{thm:CZdecRegular}. To estimate the $b_\mathrm{d}$-term, we first rewrite the $j$-th diagonal term as a sum of its restrictions to dyadic cubes in $\mathscr{D}_j$ 
$$
\bd = \sum_{j \ge 1} p_j(f-f_j)p_j= \sum_{j \ge 1} \sum_{Q \in \mathscr{D}_j} 1_Q \otimes \pi_Q(f-f_Q) \pi_Q =: \sum_{j \ge 1} \sum_{Q \in \mathscr{D}_j} b_{\mathrm{d},Q},
$$
where 
$$
f_Q = \frac{1}{|Q|} \int_Q f(x) \, dx.
$$
Observe that the above means that $b_{\mathrm{d},Q} = 1_Q \bdj $ for $Q \in \mathscr{D}_j$. Assume that we have fixed a cube $Q\in\mathscr{D}_j$. Since $\supp_{\mathbb{R}^d}(b_{\mathrm{d},Q}) \subset Q$, if $x\not\in 5Q$ and $c_Q$ denotes the center of $Q$, then we can use the kernel representation and $\mathsf{E}_j(b_{\mathrm{d},Q})= \int b_{\mathrm{d},Q} = 0$ to get
\begin{eqnarray*}
\zeta(x) T_k b_{\mathrm{d},Q}(x) \zeta(x) & = & \zeta(x) \Big( \int_{\mathbb{R}^d} k(x,y)b_{\mathrm{d},Q}(y) \, dy \Big) \zeta(x) \\
 & = & \zeta(x) \Big( \int_{Q} [k(x,y)-k(x,c_Q)] b_{\mathrm{d},Q}(y) \, dy \Big) \zeta(x).
\end{eqnarray*}
On the other hand, if $x \in 5Q$ we have $\zeta(x) T_k b_{\mathrm{d},Q}(x) \zeta(x) = 0$ by Lemma \ref{lem:zetaproydoubling}. This gives 
\begin{eqnarray*}
\varphi \big( |\zeta T_kb_{\mathrm{d},Q} \zeta| \big) \!\! & = & \!\! \varphi \big( 1_{(5Q)^c} |\zeta T(b_{\mathrm{d},Q}) \zeta| \big) \\ [5pt]
\!\! & \leq & \!\! \int_{(5Q)^c} \int_{Q} \mathrm{tr} \, \big| [k(x,y)-k(x,c_Q)] b_{\mathrm{d},Q}(y) \big| \, dy dx\\
\!\! & \leq & \!\! \int_{Q} \mathrm{tr} \, |b_{\mathrm{d},Q}(y)| \int_{(5Q)^c} \big| k(x,y)-k(x,c_Q) \big| \, dx dy \, \lesssim \, \|b_{\mathrm{d},Q}\|_1 \, = \, \|1_Q \bdj\|_1,
\end{eqnarray*}
using H\"ormander's kernel condition (weaker than our hypothesis). Finally, Lemma \ref{thm:CZdecRegular} yields
\begin{eqnarray*}
\varphi \big\{ |\zeta T_k\bd \zeta| > \lambda \big\} & \leq & \frac{1}{\lambda} \big\| \zeta T_k\bd \zeta \big\|_1 \\
& \leq & \frac{1}{\lambda} \sum_{j \ge 1} \sum_{Q \in \mathscr{D}_j} \big\| \zeta T_kb_{\mathrm{d},Q} \zeta \big\|_1 \\
& \lesssim & \frac{1}{\lambda} \sum_{j \ge 1} \sum_{Q \in \mathscr{D}_j} \|1_Q \bdj\|_1 = \frac{1}{\lambda} \sum_{j \ge 1} \|\bdj\|_1 \leq \frac{2}{\lambda} \|f\|_1.
\end{eqnarray*}

\noindent \textbf{B. Nonstandard estimates.} Recall the $j$-th term of $\boff$ $$\boffj = p_jfq_j + q_jfp_j =: \boffj^\mathrm{a} + \boffj^\mathrm{b}.$$ By symmetry, we only deal with the first term in the right hand side. As before, for each $Q \in \mathscr{D}_j$ we define the operator $\boffQ^\mathrm{a} = 1_Q \boffj^\mathrm{a} = \pi_Q f \xi_Q$ and notice that $\supp_{\mathbb{R}^d}(\boffQ^\mathrm{a}) \subset Q$ and $\int \boffQ^\mathrm{a} =0$. If $x \in 5Q$, then Lemma \ref{lem:zetaproydoubling} gives once more that $\zeta(x) T_k\boffQ^\mathrm{a}(x) \zeta(x) = 0$, while for $x \not\in 5Q$ we have
\begin{eqnarray*}
\zeta(x) T_k\boffQ^\mathrm{a}(x) \zeta(x) & = & \zeta(x) \Big( \int_{\mathbb{R}^d} k(x,y)\boffQ^\mathrm{a}(y) \, dy \Big) \zeta(x) \\
 & = & \zeta(x) \Big( \underbrace{\int_{Q} [k(x,y)-k(x,c_Q)] \pi_Q f(y) \xi_Q \, dy}_{\mathrm{B}_{Q}(x)} \Big) \zeta(x). \\
\end{eqnarray*}
If we set $K_Q(x,y) = k(x,y)-k(x,c_Q)$, this proves that 
\begin{eqnarray} \label{EqRowColumnFact} 
\hskip20pt \|\zeta T_k\boffQ^\mathrm{a} \zeta\|_1 \!\!\! & \le & \!\!\! \|1_{(5Q)^c} \mathrm{B}_Q\|_1 
\\ [5pt] \nonumber \!\!\! & \le & \!\!\! \int_{(5Q)^c} \Big\| \Big( \int_Q |K_Q(x,y)|^2 \pi_Q f(y) \pi_Q \, dy \Big)^{\frac12} \Big\|_1 dx \  \Big\| \int_Q \xi_Q f(y) \xi_Q \, dy \Big\|_{\infty}^{\frac12} 
\\ \nonumber \!\!\! & \le & \!\!\! \int_{(5Q)^c} \Big\| \Big( \int_Q |K_Q(x,y)|^2 \pi_Q f(y) \pi_Q \, dy \Big)^{\frac12} \Big\|_2 dx \  \mathrm{tr}(\pi_Q)^{\frac12} \Big\| \int_Q \xi_Q f(y) \xi_Q \, dy \Big\|_{\infty}^{\frac12}.
\end{eqnarray}
Letting $\C_{Q,\ell} = \{x: 2^\ell \ell(Q) \le |x-c_Q| \le 2^{\ell+1} \ell(Q)\}$, we may estimate 
\begin{eqnarray} \label{EqL2Horm}
\lefteqn{\hskip-10pt \int_{(5Q)^c} \Big\| \Big( \int_Q |K_Q(x,y)|^2 \pi_Q f(y) \pi_Q \, dy \Big)^{\frac12} \Big\|_2 dx} 
\\ \nonumber & \le & \sum_{\ell \ge 1} \int_{\C_{Q,\ell}} \Big( \int_Q |K_Q(x,y)|^2 \, \mathrm{tr} \big( \pi_Q f(y) \pi_Q \big) \, dy \Big)^{\frac12} dx
\\ \nonumber & \le & \sum_{\ell \ge 1} |\C_{Q,\ell}|^{\frac12} \Big( \int_{\C_{Q,\ell} \times Q} |K_Q(x,y)|^2 \, \mathrm{tr} \big( \pi_Q f(y) \pi_Q \big) \, dx dy \Big)^{\frac12}
\\ \nonumber & \le & \sum_{\ell \ge 1} \Big( \sup_{y \in Q} |\C_{Q,\ell}| \int_{\C_{Q,\ell}} |K_Q(x,y)|^2 \, dx \Big)^{\frac12} \mathrm{tr} \Big( \int_Q \pi_Q f(y) \pi_Q \, dy \Big)^{\frac12} \ \lesssim \ \varphi(1_Q p_j f)^{\frac12},
\end{eqnarray}
where the last inequality follows from the $L_2$-H\"ormander condition for the kernel given in the statement of Theorem \ref{thmDoublingCZ}. Combining \eqref{EqRowColumnFact} and \eqref{EqL2Horm}, we get the following estimate for the $\lambda$-level set of the off-diagonal terms 
\begin{eqnarray*}
\hskip26pt \varphi \big\{ |\zeta T_k \boff \zeta| > \lambda \big\} & \leq & \frac{1}{\lambda} \big\| \zeta T_k\boff \zeta \big\|_1 \\ [7pt]
& \leq & \frac{1}{\lambda} \sum_{j \ge 1} \sum_{Q \in \mathscr{D}_j} \big\| \zeta T_k (b_{\mathrm{off},Q}^{\mathrm{a}} + b_{\mathrm{off},Q}^{\mathrm{b}}) \zeta \big\|_1 \\
& \lesssim & \frac{1}{\lambda} \sum_{j \ge 1} \sum_{Q \in \mathscr{D}_j} \varphi(1_Q p_j f)^{\frac12} \mathrm{tr}(\pi_Q)^{\frac12} |Q|^{\frac12} \lambda^{\frac12} \\ & \le & \frac{1}{\lambda} \Big( \sum_{j \ge 1} \sum_{Q \in \mathscr{D}_j} \varphi(1_Q p_j f) \Big)^{\frac12} \Big( \lambda \sum_{j \ge 1} \sum_{Q \in \mathscr{D}_j} \varphi(1_Q p_j) \Big)^{\frac12} \\ & = & \frac{1}{\lambda} \varphi \big( (\1-q) f \big)^{\frac12} \big( \lambda \varphi (\1-q) \big)^{\frac12} \ \le \ \frac{1}{\lambda} \sqrt{\|f\|_1} \sqrt{\|f\|_1} \ = \ \frac{1}{\lambda} \|f\|_1. \hfill \hskip26pt \square
\end{eqnarray*}

\vskip10pt

\begin{remark}
\emph{Theorem A still holds true if we replace the Euclidean-Lebesguean space by any other metric measure space of homogeneous type. That is, the measure is doubling with respect to the given metric. To do that, it is enough to replace the dyadic filtration by the one generated by the so called Christ cubes \cite{christ1991}, an alternative detailed construction can also be found in \cite{hytonen-kairema2012}. It is easy to see that all our estimates below can be adapted to this situation just as in the scalar-valued case. Similarly, we may replace the matrix-algebra $(M_m, \mathrm{tr})$ by any other noncommutative measure space $(\M,\tau)$ or scalar-valued kernels by kernels taking values in the center of $\M$. Both of these generalizations are completely straightforward.} 
\end{remark}

\begin{remark}
\emph{Note that} \vskip-7pt $$\boff = \underbrace{\sum_{j \neq k} p_j (f - f_{j \vee k}) p_k}_{\mathrm{Easy \ term}} + \underbrace{(\1-q)fq + qf(\1-q) + \sum_{j \neq k} p_j f_{j \vee k} p_k}_{\mathrm{Former} \ g_{\mathrm{off}}}.$$ \emph{Namely, the first term above can be easily bounded for Lipschitz kernels just following the argument in \cite{cadilhac2018}. On the other hand, as it is known from \cite{parcet2009}, the second term is much harder to deal with and required pseudolocalization so far. The crucial novelty in our approach is to find a way to bound both terms together with less regularity for the kernel and a much simpler argument. The argument is quite unconventional in view of the classical CZ theory. After communicating it to Hong, it was also used in \cite{hong-lai-xu2020} to produce weak $L_1$ endpoints for maximal truncations.}
\end{remark}

\begin{remark}
{\rm
Recall that the H\"ormander integral condition on $k$ states that 
$$ \sup_{Q,y\in Q} \int_{(5Q)^c} \md{K_Q(x,y)} dx < \infty$$ in the terminology used so far.
It is the optimal condition that naturally appears in the classical scalar-valued version of Theorem A and it is still open to decide whether or not it is sufficient in our operator-valued setting. The $L_2$-condition we assume in Theorem A is stronger and becomes crucial when estimating $\|\zeta T_k\boffQ^\mathrm{a} \zeta\|_1$ in \eqref{EqL2Horm}. We found that adaptations of \eqref{EqRowColumnFact} and \eqref{EqL2Horm} assuming only $L_1$ integrability of $K_Q$ fail due to subtle noncommutative pathologies. In fact, the general strategy which consists in proving that $\sum_Q \|\zeta T_k\boffQ^\mathrm{a} \zeta\|_1 \lesssim \norm{f}_1$, 
using separately the integrability of each $K_Q$ would work (assuming the H\"ormander condition) only if $\sum_Q \| \boffQ^\mathrm{a} \|_1 \lesssim \norm{f}_1$, which cannot be expected. Hence, this question remains fully open. 
}
\end{remark}

\begin{remark}
{\rm
Another open problem related to our work above is to provide Calder\'on-Zygmund decompositions in fully noncommutative geometric settings and subsequently to establish weak type $(1,1)$ bounds for some version of singular integrals. Two possible candidates are quantum Euclidean spaces |for which good notions of distance and CZ operators were established in \cite{glez-junge-parcetMAMS}| and the hyperfinite $\mathrm{II}_1$ factor, which is naturally equipped with a dyadic filtration.   
}
\end{remark}


\subsection{Nondoubling measures} \label{section:NondoublingCZ}

Let $\mu$ be a measure on $\R^d$ of $n$-polynomial growth. In other words, if $B(x,r)$ denotes the ball centered at $x$ with radius $r$ in the Euclidean metric, then $\mu$ is a Radon measure and $0<n\leq d$ is an integer such that  
$$\mu(B(x,r)) \leq C_{\mu} r^n$$
for all radii $r>0$, $\mu$-almost every $x \in \mathrm{supp}(\mu)$ and some absolute constant $C_\mu$. On the other hand, a fixed ball $B$ is said to be $(\alpha,\beta)$-doubling if $\mu(\alpha B) \leq \beta \mu(B)$. Abundance of $(\alpha,\beta)$-doubling balls is guaranteed for appropriate values of $\alpha$ and $\beta$ under the growth condition above. In \cite{condealonso-parcet2019}, the following construction of a filtration on $\R^d$ was introduced.

\begin{theorem} \label{Weak1.dyadiccubes}
Let $\mu$ be a measure of $n$-polynomial growth on $\mathbb{R}^d$. Then, there are positive constants $\alpha,\beta > 100$ and a two-sided filtration 
$\{\Sigma_k\}_{k \in \mathbb{Z}}$ of atomic $\sigma$-algebras of $\mathrm{supp}(\mu)$ that satisfy the following properties:
\begin{itemize}
\item[\emph{i)}] The $\sigma$-algebras $\Sigma_k$ are increasingly nested.

\item[\emph{ii)}] $\bigcup_k L_\infty(\mathbb{R}^d, \Sigma_k, \mu)$ is weak-$*$ dense in $L_\infty(\mu)$.

\item[\emph{iii)}] If $Q \in \Pi(\Sigma)$, there exists an $(\alpha,\beta)$-doubling ball $B_Q$ with $B_Q \subset Q \subset 28B_Q$.

\item[\emph{iv)}] If $x \in Q \in \Pi(\Sigma)$, then $$R \ = \bigcap_{\begin{subarray}{c} S \in \Pi(\Sigma) \\ S \supsetneq Q \end{subarray}} S \qquad \Rightarrow \qquad \int_{\alpha B_R\setminus 56B_Q} \frac{d\mu(y)}{|x-y|^n}  \lesssim_{n,d,\alpha,\beta} 1.$$

\end{itemize}
\end{theorem}

Fix the measure $\mu$, and the values $\alpha$ and $\beta$ guaranteed by Theorem \ref{Weak1.dyadiccubes}. On $\A = L_\infty(\mu) \otimes M_m$, we consider the filtration $\{\A_j\}_{j\in\Z}$ generated by the filtration $\{\Sigma_j\}_{j\in\Z}$, as explained at the beginning of this section. Consider a CZ operator $T_k$ with representation
$$
T_k f(x) = \int_{\R^d} k(x,y) f(y) d\mu(y) \quad \mbox{for} \quad x \notin \supp_{\R^d}f.
$$
As usual, we fix a compactly supported $f \in L_1(\A)_+$ and fix $\lambda >0$. Using Cuculescu's projections for $(f,\lambda)$ we need once more an auxiliary projection, similar to the one that we used in the doubling case. We define
$$
\zeta := \1_\A - \bigvee_{j \ge 1} \bigvee_{Q \in  \Pi(\Sigma_j)} 1_{\alpha B_Q} \pi_Q.
$$
Since $B_Q$ is $(\alpha,\beta)$-doubling for all $Q\in \Pi(\Sigma)$, we have a result similar to Lemma \ref{lem:zetaproydoubling} and whose proof follows once again as in \cite[Lemma 3.4]{cadilhac2018}. More precisely, the above defined projection $\zeta$ satisfies the inequality  
$$
\varphi(\mathbf{1}_\A - \zeta) \lesssim \frac{\|f\|_{L_1(\A)}}{\lambda}.
$$
Moreover, we have $\zeta(x)p_j(y) = p_j(y)\zeta(x) = 0$ whenever $y \in Q \in \Pi(\Sigma_j)$ and $x \in \alpha B_Q$.


\demB Given $(f,\lambda)$ as above and according to the properties of the projection $\zeta$, we are reduced as in Theorem \ref{thmDoublingCZ} to prove the following inequality for the different parts of the nondoubling CZ decomposition of $f=g+\bd+\boff$ at level $\lambda$
$$
\varphi \big\{ |T_kg|>\lambda \big\} + \varphi \big\{|\zeta T_k\bd \zeta| > \lambda \big\} + 
\varphi \big\{ |\zeta T_k\boff \zeta| > \lambda \big\} \, \lesssim \, \frac{\|f\|_{L_1(\A)}}{\lambda}.
$$ 
According to the inequality $\|g\|_2^2 \lesssim \lambda \|f\|_1$ from Lemma \ref{thm:CZdecNonRegular}, the first term above is estimated just as we did in Theorem \ref{thmDoublingCZ}. The term $\varphi (\{ |\zeta T\boff \zeta| > \lambda\})$ is also estimated as in the proof of Theorem \ref{thmDoublingCZ} just replacing $5Q$ by $\alpha B_Q$ and using the coronas $\C_{Q,\ell} = \{x: 2^\ell r(B_Q) \le |x-c_Q| \le 2^{\ell+1} r(B_Q)\}$ defined in the statement of Theorem \ref{thmNonDoublingCZ}. To estimate $\varphi\left(\zeta |T\bd| \zeta > \lambda \right)$, we write
$$
\bd = \sum_{j \ge 1} p_j f p_j - \mathsf{E}_{j-1}(p_jfp_j) = \sum_{j \ge 1} \sum_{Q \in \Pi(\Sigma_j)} \Big( \pi_Q f \pi_Q 1_Q  - \frac{\mu(Q)}{\mu(\widehat{Q})} \pi_Q f_Q \pi_Q  1_{\widehat{Q}} \Big) =: \sum_{j \ge 1} \sum_{Q \in \Pi(\Sigma_j)}  b_Q. 
$$
Above, $\widehat{Q}$ denotes the unique atom $R \in \Pi(\Sigma_{j-1})$ containing $Q$. If the filtration were the dyadic one, it would be that dyadic parent of $Q$. Let us work with one fixed $Q\in \Pi(\Sigma_j)$. Notice that $\supp (b_Q) \leq 1_{\widehat{Q}} \pi_Q$. We consider different situations:
\begin{itemize}
\item If $x\not\in \alpha B_{\widehat{Q}}$ we use the mean $0$ of $b_Q$ to get 
\begin{eqnarray*}
\zeta(x) T_kb_Q(x) \zeta(x) & = & \zeta(x) \int_{\Omega} k(x,y)b_Q(y) \, d\mu(y) \zeta(x). \\
 & = & \zeta(x) \int_{\widehat{Q}} [\underbrace{k(x,y)-k(x,x_Q)}_{K_Q(x,y)}] b_Q(y) \, d\mu(y) \zeta(x). \\
\end{eqnarray*}
By Fubini's theorem and H\"ormander's regularity (weaker than $L_2$-H\"ormander), we obtain
\begin{eqnarray*}
\hskip30pt \int_{(\alpha B_{\widehat{Q}})^c} \mathrm{tr} |\zeta T_kb_Q \zeta| \, d\mu \!\!\! & \lesssim & \!\!\! \int_{(\alpha B_{\widehat{Q}})^c} \mathrm{tr} \left| \int_{\widehat{Q}} K_Q(x,y) b_Q(y) d\mu(y) \right| \, d\mu(x) \\
& \leq & \sup_{y \in \widehat{Q}} \int_{(\alpha B_{\widehat{Q}})^c} \big| k(x,y) - k(x,c_Q) \big| \, dx \int_{\widehat{Q}}  \mathrm{tr} | b_Q(y) | \, dy \, \lesssim \, \|b_Q\|_{L_1(\A)}.
\end{eqnarray*}
\item If $x \in \alpha B_{\widehat{Q}} \setminus \alpha B_{Q}$ we use the size kernel condition in the statement to compute
\begin{eqnarray*}
\hskip33pt \int_{\alpha B_{\widehat{Q}}\setminus \alpha B_{Q}} \mathrm{tr} |\zeta T_k(b_Q) \zeta| \, d\mu \!\!\! & \leq & \!\!\!  \int_{\alpha B_{\widehat{Q}}\setminus \alpha B_{Q}} \int_{Q} \frac{1}{|x-y|^n} \mathrm{tr} \big( \pi_Q f(y) \pi_Q \big) \, d\mu(y)d\mu(x) \\
\!\!\! & + & \!\!\! \int_{\alpha B_{\widehat{Q}}} \mathrm{tr} \Big| T_k \Big(\frac{1_{\widehat{Q}}}{\mu(\widehat{Q})} \pi_Q \int_Qf d\mu \ \pi_Q \Big)(x) \Big| \, d\mu(x) \, =: \, \mathrm{A}_Q + \mathrm{B}_Q.
\end{eqnarray*}
Theorem \ref{Weak1.dyadiccubes} iv) gives $\mathrm{A}_Q \lesssim \varphi(1_Q \pi_Q f)$. On the other hand
 \begin{eqnarray*}
 \mathrm{B}_Q & = &\frac{1}{\mu(\widehat{Q})} \int_{\alpha B_{\widehat{Q}}} \big| T_k (1_{\widehat{Q}}) \big| \, d\mu \ \mathrm{tr} \Big( \pi_Q \int_Qf d\mu \ \pi_Q \Big) \\
 & \leq & \frac{1}{\mu(\widehat{Q})} \mu(\alpha B_{\widehat{Q}})^{\frac12} \Big( \int_{\alpha B_{\widehat{Q}}} \big| T_k (1_{\widehat{Q}}) \big|^2 \, d\mu \Big)^{\frac12} \varphi(1_Q p_j f) \\ [3pt]
 & \lesssim & \frac{1}{\mu(\widehat{Q})} \mu(\alpha B_{\widehat{Q}})^{\frac12} \mu(\widehat{Q})^{\frac12} \varphi(1_Q p_j f) \ \, \le \, \ \beta \varphi(1_Q p_j f). 
\end{eqnarray*}

\item Finally, if $x \in \alpha B_{Q}$, we know from the properties of $\zeta$ that $\zeta(x) T(b_Q)(x) \zeta(x) = 0$.
\end{itemize}
Finally, we add everything up and use the $L_1$ estimate of $\bd$:
\begin{eqnarray*}
\hskip38pt \varphi \big\{ \zeta |T_k\bd| \zeta > \lambda \big\} & \leq & \frac{1}{\lambda}\left\| \zeta |T_k\bd|\zeta \right\|_{L_1(\A)} \\ [7pt]
& \leq &  \frac{1}{\lambda} \sum_{j \ge 1} \sum_{Q \in \Pi(\Sigma_j)} \big\| \zeta |T_kb_Q|\zeta \big\|_{L_1(\A)} \\
& \lesssim &  \frac{1}{\lambda} \sum_{j \ge 1} \sum_{Q \in \Pi(\Sigma_j)} \big[ \|1_Q \bdj\|_{L_1(\A)} + \varphi(1_Q p_j f) \big] \ \lesssim \ \|f\|_{L_1(\A)}. \hfill \hskip39pt \square
\end{eqnarray*}

\vskip10pt

\begin{remark}
\emph{Using the $n$-polynomial growth of $\mu$, it is easily checked that our kernel regularity assumption is weaker than Lipschitz regularity, which for this class of measures is formulated as $$\big| k(x,y) - k(x,z) \big| \le \frac{|y-z|^\gamma}{|x-y|^{n+\gamma}} \quad \mbox{for some} \quad \gamma > 0.$$}
\end{remark}

\begin{remark}
\emph{The generalizations of Theorem \ref{thmDoublingCZ} concerning other noncommutative measure spaces $(\M,\tau)$ or kernels taking values in the center $\mathcal{Z}(\M)$ still apply for Theorem \ref{thmNonDoublingCZ}. In the latter case, we need to use from \cite[Remark 2.4]{junge-mei-parcetSmooth} that the $L_2$ boundedness of $T_k$ implies for this class of kernels that it maps $L_\infty(\M;L_2^c(\mu))$ onto itself. In other words} 
$$
\Big\|\int_{\R^d} |T_kf|^2 \, d\mu \Big\|_{\M}^{\frac{1}{2}} \lesssim \Big\| \int_{\R^d} |f|^2 \, d\mu \Big\|_{\M}^{\frac{1}{2}}.
$$
\emph{It is unknown whether Theorem \ref{thmNonDoublingCZ} holds for noneuclidean measure spaces of $n$-polynomial growth.}
\end{remark}


\section{\bf Fourier and Schur multipliers} \label{sectionFourier}

Let $\G$ be a locally compact group. A length $\psi:  \G \to \R_+$ is any continuous function satisfying $\psi(e) = 0$ and $\psi(g) = \psi(g^{-1})$, where $e$ denotes the unit in $\G$ and $g$ is a generic element. It is called conditionally negative when $\sum a_g \overline{a_h} \psi(gh^{-1}) \le 0$ for finite families of coefficients $a_g$ satisfying $\sum a_g = 0$. An  orthogonal cocycle is given by a real Hilbert space $\mathcal{H}$, an orthogonal representation $\alpha: \G \to O(\mathcal{H})$ and a map $\beta: \G \to \mathcal{H}$ satisfying $\alpha_g(\beta(h)) = \beta(gh) - \beta(g)$. Conditionally negative lengths are in one-to-one correspondence with orthogonal cocycles by $\psi(g) = |\beta(g)|^2$. The dimension of $\psi$ is thus defined as $\dim \mathcal{H}$ for the only cocycle related to $\psi$. A more detailed presentation and references to all these concepts and relations may be found in \cite{junge-mei-parcetSmooth} and \cite{junge-mei-parcetJEMS}.

A crucial advantage of Theorem \ref{thmDoublingCZ} in comparison with its ancestors in \cite{cadilhac2018,parcet2009} is that the kernel regularity condition now includes H\"ormander-Mikhlin multipliers. This assertion is readily implied by \cite[Lemma 1]{Wheeden}. We give a precise statement for future reference. 

\begin{lemma} \label{LemmaHMCZ}
If $M: \R^n \to \mathbb{C}$ satisfies the H\"ormander-Mikhlin condition 
$$
\big| \partial_\xi^\gamma M(\xi) \big| \, \lesssim \, |\xi|^{-|\gamma|} \quad \mbox{for all} \quad 0 \le |\gamma| \le \Big[ \frac{n}{2} \Big] + 1,
$$
then the kernel $k(x,y) = \widehat{M}(x-y)$ satisfies the $L_2$-H\"ormander condition in Theorem \ref{thmDoublingCZ} for $d=n$. 
In particular, for any noncommutative measure space $\mathcal{M}$
$$
\norm{T_M \otimes \mathrm{Id}}_{L_1(L_\infty(\R^n)\overline{\otimes}\mathcal{M}) \to L_{1,\infty}(L_\infty(\R^n)\overline{\otimes}\mathcal{M})} \, \lesssim \ \sup_{\xi \neq 0} \sum_{\md{\gamma} \leq [n/2] + 1} \md{\xi}^{\md{\gamma}}\big| \partial_\xi^\gamma M(\xi) \big| =:  \norm{M}_{\mathrm{HM}}.
$$

\end{lemma}

\subsection{Proof of Theorem \ref{thmFourierMult} i)}


In this paragraph, we will prove our statement for Schur multipliers of Mikhlin type. More precisely, we shall establish weak $L_1$ bounds for Schur multipliers $S_M$ with symbol $(g,h) \mapsto M(gh^{-1}) = m \circ \psi (gh^{-1})$ for some conditionally negative length $\psi: \G \to \R_+$ and a spectral multiplier $m: \R_+ \to \mathbb{C}$ satisfying the Mikhlin condition for $\dim \psi = n$. Let assume that $\G$ is discrete for simplicity, although (as it will be remarked below) our argument generalizes to every locally compact group as long as it satisfies the hypotheses in Theorem \ref{thmFourierMult}. In the line of \cite{CasdlS,parcet-ricard-salle} we define $\pi: S_p(\G) \to L_\infty(\R^n; S_p(\G))$ by $$\pi(A) =  \Big( \exp_{\alpha_{g^{-1}}(\beta(gh^{-1}))} A_{gh} \Big)_{g,h} = \Big( \exp_{-\beta(g^{-1})} A_{gh} \exp_{\beta(h^{-1})} \Big)_{g,h} = u^* (\1 \otimes A) u$$ with $\exp_\xi(x) = e^{2 \pi i \langle x,\xi \rangle}$ and $u = \mathrm{diag} (\exp_{\beta(g^{-1})})$. Next, we claim that 
\begin{eqnarray} \label{EqClaim1}
\|S_M(A)\|_{S_{1,\infty}(\G)} & \lesssim & \lim_{\varepsilon \to 0} \big\| T_{\widetilde{m}}(\gamma_\varepsilon \pi(A)) \big\|_{L_{1,\infty}(\A)} \\ \nonumber & \lesssim & \lim_{\varepsilon \to 0} \big\| \gamma_\varepsilon \pi(A) \big\|_{L_1(\A)} \ = \ \|A\|_{S_1(\G)}.
\end{eqnarray}
for the Fourier multiplier associated to the lifted symbol $\widetilde{m}: \R^n \to \mathbb{C}$ determined by $\widetilde{m}(\xi) = m(|\xi|^2)$ and the family of $L_1$-normalized gaussians $\gamma_\varepsilon(x) = ( \varepsilon/\pi )^{\frac{n}{2}} \exp(- \varepsilon |x|^2)$. This implies the statement in the Banach space setting, but the same argument applies after matrix amplification. Therefore it suffices to justify claim \eqref{EqClaim1} to complete the proof of Theorem \ref{thmFourierMult} i). 

It is clear that $\| \gamma_\varepsilon \pi(A) \|_{L_1(\A)} = \| \gamma_\varepsilon \otimes A \|_{L_1(\A)} = \|A\|_{S_1(\G)}$ for all $\varepsilon > 0$. This justifies the last identity in \eqref{EqClaim1}. On the other hand, it is straightforward to show that the hypotheses in Theorem \ref{thmFourierMult} imply that the lifted symbol $\widetilde{m}$ satisfies as well the Mikhlin condition in $n$ variables. More precisely
$$
\big| \partial_\xi^\gamma \widetilde{m}(\xi) \big| \, \lesssim \, |\xi|^{-|\gamma|} \quad \mbox{for all} \quad 0 \le |\gamma| \le \Big[ \frac{n}{2} \Big] + 1.
$$
In addition, it follows from Lemma \ref{LemmaHMCZ} that the Fourier multiplier $T_{\widetilde{m}}$ is a CZ operator associated to a kernel $k$ satisfying the $L_2$-H\"ormander condition from Theorem \ref{thmDoublingCZ}. Thus, the second inequality in claim \eqref{EqClaim1} is a direct consequence of Theorem \ref{thmDoublingCZ}. Note that Lipschitz regularity from \cite{cadilhac2018,parcet2009} is not enough at this point. It remains to justify the first inequality in claim \eqref{EqClaim1}. To that end we consider the  gaussian probability measures $d\sigma_\varepsilon(x) = \gamma_\varepsilon(x) dx$ and the algebra of matrix-valued functions $\A_\varepsilon = L_\infty(\R^n, \sigma_\varepsilon) \bar\otimes \B(\ell_2(\G))$. We shall prove the remaining inequality in several steps as follows 
\begin{eqnarray} \label{EqClaim2}
\|S_M(A)\|_{S_{1,\infty}(\G)} & =_\mathbf{a} & \lim_{\varepsilon \to 0} \big\| \pi(S_M(A)) \big\|_{L_{1,\infty}(\A_\varepsilon)}
\\ \nonumber & =_\mathbf{b} & \lim_{\varepsilon \to 0} \big\| \hskip1pt T_{\widetilde{m}}(\pi(A)) \hskip0.5pt \big\|_{L_{1,\infty}(\A_\varepsilon)} 
\\ \nonumber & \lesssim_\mathbf{c} & \lim_{\varepsilon \to 0} \big\| \gamma_\varepsilon T_{\widetilde{m}}(\pi(A)) \big\|_{L_{1,\infty}(\A)} 
\\ \nonumber & \lesssim_\mathbf{d} & \lim_{\varepsilon \to 0} \big\| T_{\widetilde{m}}(\gamma_\varepsilon \pi(A)) \big\|_{L_{1,\infty}(\A)}. 
\end{eqnarray}   

\noindent \textbf{Proof of \eqref{EqClaim2}$_{\mathbf{a}}$.} Just note that 
\begin{eqnarray*}
\|S_M(A)\|_{S_{1,\infty}(\G)} & = & \sup_{\lambda > 0} \, \lambda \int_{\R^n} \mathrm{tr} \big\{ |\1 \otimes S_M(A)| > \lambda \big\} \, d\sigma_\varepsilon
\\ & = & \sup_{\lambda > 0} \, \lambda \int_{\R^n} \mathrm{tr} \big\{ |u^* (\1 \otimes S_M(A)) u| > \lambda \big\} \, d\sigma_\varepsilon \ = \ \big\| \pi(S_M(A)) \big\|_{L_{1,\infty}(\A_\varepsilon)}.
\end{eqnarray*}

\noindent \textbf{Proof of \eqref{EqClaim2}$_{\mathbf{b}}$.} Just note that 
\begin{eqnarray*}
T_{\widetilde{m}}(\pi(A)) & = & \Big( T_{\widetilde{m}}(\exp_{\alpha_{g^{-1}}(\beta(gh^{-1}))} A_{gh}) \Big)_{g,h}
\\ & = & \Big( \widetilde{m}(\alpha_{g^{-1}}(\beta(gh^{-1}))) \exp_{\alpha_{g^{-1}}(\beta(gh^{-1}))} A_{gh} \Big)_{g,h}
\\ & = & \Big( m \circ \psi (gh^{-1}) \exp_{\alpha_{g^{-1}}(\beta(gh^{-1}))} A_{gh} \Big)_{g,h} \ = \ \pi(S_M(A)).
\end{eqnarray*}

\noindent \textbf{Proof of \eqref{EqClaim2}$_{\mathbf{c}}$.} By the above identity, it suffices to prove 
$$
\| \pi(B) \|_{L_{1,\infty}(\A_\varepsilon)} \, \lesssim \, \| \gamma_\varepsilon \pi(B) \|_{L_{1,\infty}(\A)} \quad \mbox{for all} \quad \varepsilon > 0.
$$
To that end, let us introduce two parameters defined as follows:
\begin{itemize}
\item Pick $\lambda_B > 0$ such that $\|B\|_{S_{1,\infty}(\G)} \le 2 \lambda_B \mathrm{tr} \big\{ |B| > \lambda_B \big\}$.

\item Pick $R_\varepsilon = \sqrt{\frac{\log 2}{\varepsilon}}$, so that $\gamma_\varepsilon (x) = \frac12 \gamma_\varepsilon(0)$ for all $x$ with $|x| = R_\varepsilon$.
\end{itemize}
Then, we are in position to estimate the weak $L_1$ norm of $\gamma_\varepsilon \pi(B)$ from below 
\begin{eqnarray*}
\big\| \gamma_\varepsilon \pi(B) \big\|_{L_{1,\infty}(\A)} & = & \sup_{\lambda > 0} \, \lambda \int_{\R^n} \mathrm{tr} \big\{ |\gamma_\varepsilon u^*(\1 \otimes B)u| > \lambda \big\} \, dx \\ 
& \ge & \int_{B(0,R_\varepsilon)} \frac{\lambda_B \gamma_\varepsilon(0)}{2 \gamma_\varepsilon(x)} \mathrm{tr} \Big\{ |B| > \frac{\lambda_B \gamma_\varepsilon(0)}{2 \gamma_\varepsilon(x)} \Big\} \, d\sigma_\varepsilon(x) \\ 
& \ge & \int_{B(0,R_\varepsilon)} \frac{\lambda_B}{2} \mathrm{tr} \Big\{ |B| > \lambda_B \Big\} \, d\sigma_\varepsilon(x) \ \ge \ \frac14 \sigma_\varepsilon(B(0,R_\varepsilon)) \|B\|_{S_{1,\infty}(\G)}.
\end{eqnarray*}
However, it is not difficult to calculate that $$\sigma_\varepsilon(B(0,R_\varepsilon)) = \frac{|\mathbb{S}^{n-1}|}{(4\pi)^{n/2}} \quad \mbox{for all} \quad \varepsilon > 0.$$ Therefore, the desired inequality follows since $\| \pi(B) \|_{L_{1,\infty}(\A_\varepsilon)} = \|B\|_{S_{1,\infty}(\G)}$ for all $\varepsilon > 0$. 

\vskip5pt

\noindent \textbf{Proof of \eqref{EqClaim2}$_{\mathbf{d}}$.} We will prove the identity 
\begin{equation} \label{EqClaim3}
\lim_{\varepsilon \to 0} \big\| T_{\widetilde{m}}(\gamma_\varepsilon \pi(A)) - \gamma_\varepsilon T_{\widetilde{m}}(\pi(A)) \big\|_{L_{1,\infty}(\A)} = 0.
\end{equation} 
This readily implies claim \eqref{EqClaim2} by the quasi-triangular inequality in weak $L_1$ spaces. By density in $S_1(\G)$, we may assume from the beginning that $A_{gh} \neq 0$ for finitely many entries $(g,h)$. Again by the quasi-triangle inequality, it suffices to prove \eqref{EqClaim3} entrywise. In other words, if we set $\xi_{gh} = \beta(h^{-1}) - \beta(g^{-1})$, then we must prove that 
$$
\lim_{\varepsilon \to 0} \big\| \big( \underbrace{T_{\widetilde{m}}(\gamma_\varepsilon \exp_{\xi_{gh}}) - m \circ \psi(gh^{-1}) \gamma_\varepsilon \exp_{\xi_{gh}}}_{\Sigma_{gh}^\varepsilon} \big) \otimes e_{gh} \big\|_{1,\infty} = \lim_{\varepsilon \to 0} \|\Sigma_{gh}^\varepsilon \|_{1,\infty} = 0 \quad \mbox{for all $g,h \in \G$}.
$$
If we set $\Phi(\xi) = \widetilde{m}(\xi + \xi_{gh})$, we may easily construct $\Psi_\varepsilon: \R^n \to \R_+$ with:
\begin{itemize}
\item $\Psi_\varepsilon$ smooth, 

\vskip2pt

\item $\Psi_\varepsilon(\xi)=1$ when $|\xi| < \varepsilon^{\frac14}$ and $\Psi_\varepsilon(\xi) = 0$ when $|\xi| > 2 \varepsilon^{\frac14}$,


\item $\displaystyle \lim_{\varepsilon \to 0} \, \sup_{\xi \neq 0} \sum_{|\gamma| \le [n/2]+1} \md{\xi}^{\md{\gamma}}\big| \partial_\xi^\gamma \big[ \Psi_\varepsilon(\xi) (\Phi(\xi) - \Phi(0)) \big] \big| = \lim_{\varepsilon \to 0} \big\| \Psi_\varepsilon (\Phi - \Phi(0)) \big\|_{\mathrm{HM}}= 0$.
\end{itemize}
Then, using the identity $T_M(f \exp_{\xi_0}) = T_{M(\cdot + \xi_0)}(f) \exp_{\xi_0}$, we may continue our estimate as follows 
\begin{eqnarray*}
\hskip20pt \lim_{\varepsilon \to 0} \|\Sigma_{gh}^\varepsilon \|_{1,\infty} & = & \lim_{\varepsilon \to 0} \big\| T_\Phi(\gamma_\varepsilon) - \Phi(0) \gamma_\varepsilon \big\|_{1,\infty} \\ \nonumber & = & \lim_{\varepsilon \to 0} \big\| T_\Phi \big( \gamma_\varepsilon -  T_{\Psi_\varepsilon} (\gamma_\varepsilon) \big) + T_{\Psi_\varepsilon (\Phi-\Phi(0))} (\gamma_\varepsilon) + \Phi(0) \big( T_{\Psi_\varepsilon} (\gamma_\varepsilon) - \gamma_\varepsilon \big) \big\|_{1,\infty} \\ \nonumber & \lesssim & \lim_{\varepsilon \to 0} \|\Phi\|_{\mathrm{HM}} \big\| \gamma_\varepsilon -  T_{\Psi_\varepsilon} (\gamma_\varepsilon) \big\|_1 + \lim_{\varepsilon \to 0} \big\| \Psi_\varepsilon (\Phi-\Phi(0)) \big\|_{\mathrm{HM}} \|\gamma_\varepsilon\|_1.
\end{eqnarray*}
Since $\|\gamma_\varepsilon\|_1 = 1$, the last limit on the right hand side vanishes by construction of $\Psi_\varepsilon$. On the other hand, we have $\|\Phi\|_{\mathrm{HM}} \lesssim \|\widetilde{m}\|_{\mathrm{HM}} \lesssim \|m\|_{\mathrm{HM}} < \infty$ by hypothesis. Therefore, it remains to prove that $\| \gamma_\varepsilon -  T_{\Psi_\varepsilon} (\gamma_\varepsilon) \|_1$ gets arbitrarily small as $\varepsilon \to 0$. Since $\widehat{\gamma_\varepsilon}(\xi) = \exp(-\pi|\xi|^2/\varepsilon)$, we may write this $L_1$-norm as follows
\begin{eqnarray*}
\big\| \gamma_\varepsilon -  T_{\Psi_\varepsilon} (\gamma_\varepsilon) \big\|_1 & = & \int_{\R^n} \Big| \int_{\R^n} (1 - \Psi_\varepsilon(\xi)) e^{-\pi |\xi|^2/\varepsilon} e^{2\pi i \langle x,\xi \rangle} \, d\xi \Big| \, dx
\\ & = & \int_{\R^n} \Big| \int_{\R^n} \underbrace{(1 - \Psi_\varepsilon(\varepsilon^{\frac12} \eta)) e^{-\pi |\eta|^2}}_{\Lambda_\varepsilon(\eta)} e^{2\pi i \langle y,\eta \rangle} \, d\eta \Big| \, dy \ = \ \| \widehat{\Lambda}_\varepsilon \|_1.
\end{eqnarray*}
We know that $\Lambda_\varepsilon \to 0$ pointwise as $\varepsilon \to 0$. By the dominated convergence theorem, this also holds in $L_1$-norm. According to Hausdorff-Young inequality, this implies that its Fourier transform converges to $0$ in the $L_\infty$-norm. We also get 
$$
\big\| |x|^{|\gamma|} \widehat{\Lambda}_\varepsilon \big\|_\infty \, \asymp \, \big\| \widehat{\partial_\xi^\gamma \Lambda}_\varepsilon \big\|_\infty \, \le \, \big\| \partial_\xi^\gamma \Lambda_\varepsilon \big\|_1
$$
which converges to $0$ once again by dominated convergence. Taking $|\gamma| = n+1$, we conclude. \fin  
 
\begin{remark}
{\emph{The above argument can be easily modified to hold for arbitrary locally compact groups. Indeed, the argument follows verbatim up to inequality \eqref{EqClaim2}$_{\mathbf{d}}$. At this point, we may assume by density that $A$ has a kernel in $\C_c(\G \times \G)$ and argue as in \cite[Theorem 1.19]{LafforgueSalle} to approximate $\pi(A)$ in $L_1(\B(L_2(\G,\mu))$ by $$\sum_{j,k} \exp_{\xi_{g_jh_k}} A_{g_jg_k} e_{g_jg_k}$$ 
in the space $L_1(\B(L_2(\G,\mu'))$, where $\mu'$ is certain finitely supported measure on $\G$. Then, the argument follows by estimating the weak $L_1$ norm of the functions $\Sigma_{g_jh_k}^\varepsilon$ exactly as done above.}}
\end{remark}

\begin{remark}
\emph{The above argument can also be modified to give a weak $L_1$ form of de Leeuw's compactification theorem \cite{dL}. More precisely, let $\R_{\mathrm{bohr}}$ denote Bohr's compactification of the real line $\R$ and consider a symbol $M: \R \to \mathbb{C}$ giving rise to a weak-$L_1$-bounded Fourier multiplier $T_M: L_1(\R) \to L_{1,\infty}(\R)$. We may also regard $M$ as a symbol on the discrete real line $\R_{\mathrm{disc}}$, the Pontryagin predual of $\R_{\mathrm{bohr}}$. Then, we can use the above ideas to prove}
$$
\big\| T_M: L_1(\R_{\mathrm{bohr}}) \to L_{1,\infty}(\R_{\mathrm{bohr}}) \big\| \, \le \, \big\| T_M: L_1(\R) \to L_{1,\infty}(\R) \big\|.
$$
\emph{In fact, arguing as in \cite[Section 2.5]{ParRog} the same holds for group algebras over unimodular groups.}
\end{remark}

\subsection{Proof of Theorem \ref{thmFourierMult} ii)}

The discreteness assumption in Theorem \ref{thmFourierMult} ii) is to avoid further considerations that would lead us too far from our central topic.
On the other hand, amenability is strongly linked to the transference methods that we need. Almost every form of transference since the pioneer work of Cotlar or Calder\'on involve some kind of amenability. In this line, we shall generalize the tight connection between Fourier and Schur $L_p$-multipliers \cite{CasdlS, NeuRic} through the key inequality for $\G$ amenable 
$$\|T_M\|_{\mathrm{cb}(L_p(\V))} \le \|S_M\|_{\mathrm{cb}(S_p(\G))} \quad \mbox{when} \quad 1 \le p \le \infty.$$ In conjunction with Theorem \ref{thmFourierMult} i), the result below (of independent interest) yields Theorem \ref{thmFourierMult} ii).

\begin{theorem}\label{thmTransfer}
Let $\G$ be a  discrete group, then 
\begin{eqnarray*}
\lefteqn{\hskip-40pt \Big\| \mathrm{id} \otimes T_M: L_1(M_m \otimes \V) \to L_{1,\infty}(M_m \otimes \V) \Big\|} \\ \hskip40pt & \le & \, \Big\| \mathrm{id} \otimes S_M: S_1(\ell_2(m) \otimes \ell_2(\G)) \to S_{1,\infty}(\ell_2(m) \otimes \ell_2(\G)) \Big\|.
\end{eqnarray*} 
\end{theorem}

\dem The argument follows very much the ideas in \cite{CasdlS, NeuRic, parcet-ricard-salle}. By amenability, we may consider a F\o lner net $(\Lambda_\alpha)_\alpha$. In other words, a family of finite non empty sets satisfying that $\lim_\alpha \md{\Lambda_\alpha \cap g \Lambda_\alpha} / \md{\Lambda_\alpha} = 1$ for all $g \in \G$. Consider the canonical inclusion map $j: \V\to \B(\ell_2(\G))$ given by 
$$
j(f) = \sum_{g,h\in G} \widehat{f}(gh^{-1})e_{gh}.
$$
The map $j$ is a $*$-homomorphism. It fails to be $L_p$-bounded, but
$$
\|f\|_{L_p(\V)} \, = \, \lim_{\alpha} \|j_{\alpha}(f)\|_{S_p(\Lambda_\alpha)}
$$
$$
\mbox{where} \quad j_{\alpha}: \V \ni f \mapsto \pi_\alpha \hskip2pt j(f) \hskip1pt \pi_\alpha \in S_\infty(\Lambda_\alpha),
$$ \vskip5pt
\noindent $S_\infty(\Lambda_\alpha)$ is equipped with its normalised trace $\mathrm{tr}/\md{\Lambda_\alpha}$ and $\pi_\alpha$ denotes the orthogonal projection onto the subspace $\ell_2(\Lambda_\alpha)$. Consider a tracial ultraproduct $\big( \prod_{\mathcal{U}} S_{\infty}(\Lambda_\alpha),\tau_{\mathcal{U}} \big)$ of the Schatten classes $S_\infty(\Lambda_\alpha)$. The map
\begin{align*}
J \colon \V \ni f & \mapsto \overline{(j_\alpha(f))}_\alpha \in \prod_{\mathcal{U}} S_{\infty}(\Lambda_\alpha),
\end{align*}
is a trace preserving completely isometric embedding. It is clearly completely positive and unital so we should only check that it is multiplicative, which essentially comes down to proving that for any $f_1,f_2 \in \V$, we have $\tau(f_1f_2) = \tau_{\mathcal{U}}(J(f_1)J(f_2))$. It suffices to check this equality for elements $f_1 = \lambda(g_1),f_2 = \lambda(g_2)$, $g_1,g_2 \in\G$. If $g_1 \neq g_2^{-1}$, both terms are $0$ and otherwise, for any $g\in\G$
\begin{eqnarray*}
\tau_{\mathcal{U}} \big( J(\lambda(g)) J(\lambda(g^{-1})) \big) & = &
\lim_\alpha \dfrac{1}{\md{\Lambda_\alpha}} \mathrm{tr} \big( \pi_\alpha j( \lambda(g)) \pi_\alpha j(\lambda(g^{-1})) \pi_\alpha \big) \\
& = & \lim_\alpha \dfrac{1}{\md{\Lambda_\alpha}} \sum_{h,h'\in \Lambda_{\alpha}} \delta_{h',gh} \ = \ \lim_\alpha \dfrac{\md{\Lambda_\alpha \cap g\Lambda_\alpha}}{\md{\Lambda_\alpha}} \ = \ 1.
\end{eqnarray*}
To conclude, observe that $J$ intertwines $T_M$ and $\widetilde{S}_M := \Pi_\mathcal{U} S_M$. Hence for any $f\in L_1(\V)$
\begin{eqnarray*}
\norm{T_M(f)}_{1,\infty} & = & \big\| J(T_M(f)) \big\|_{1,\infty} \ = \ \big\| \widetilde{S}_M(J(f)) \big\|_{1,\infty} \\
& \leq & \big\| \widetilde{S}_M: L_1 \to L_{1,\infty} \big\| \, \|J(f)\|_1 \ \leq \ \big\| S_M: L_1 \to L_{1,\infty} \big\| \, \|f\|_1.
\end{eqnarray*}
\fin

\begin{remark}
{\rm 
We expect Theorem \ref{thmTransfer} above to hold for any locally compact unimodular amenable group. However, to the best of our knowledge, this cannot be easily infered from the literature and only the discrete case lies within the scope of this paper. Note also that alternative arguments in the spirit of \cite[Theorem 2.1]{NeuRic} could have been provided. 
}
\end{remark}

\begin{remark}
\emph{It is quite standard to use our results so far to establish a $L_1$ endpoint for the Littlewood-Paley theorem in group von Neumann algebras \cite{junge-mei-parcetSmooth}. More precisely, if $(m_j)_{j \in \Z}$ is a Littlewood-Paley partition of unity in $\R_+$ and $\psi: \G \to \R_+$ is a conditionally negative length, then we set $M_j = m_j \circ \psi$ and obtain}
\begin{itemize}
\item \emph{LP for Schur multipliers} 
$$
\hskip30pt \inf_{S_{M_j}A = R_jA + L_jA} \Big\| \Big( \sum_{j \in \Z} R_jAR_jA^* + L_jA^*L_jA \Big)^{\frac12} \Big\|_{S_{1,\infty}(\G)} \, \le \, \Big( \sum_{j \in \Z} \|m_j\|_{\mathrm{HM}}^2 \Big)^{\frac12} \|A\|_{S_1(\G)}.
$$   

\item \emph{LP for Fourier multipliers}
$$
\hskip30pt \inf_{T_{M_j}f = A_jf + B_jf} \Big\| \Big( \sum_{j \in \Z} A_jfA_jf^* + B_jf^*B_jf \Big)^{\frac12} \Big\|_{L_{1,\infty}(\V)} \, \le \, \Big( \sum_{j \in \Z} \|m_j\|_{\mathrm{HM}}^2 \Big)^{\frac12} \|f\|_{L_1(\V)}.
$$   
\end{itemize}
\emph{As usual, these results admit matrix amplifications and hold for l.c. (respectively amenable) groups.}
\end{remark}

\noindent {\bf Acknowledgements.} J. Parcet and J.M. Conde-Alonso were partially supported by Spanish Grant PID2019-107914GB-I00 and Severo Ochoa Programme for Centres of Excellence in R\&D CEX2019-000904-S. J.M. Conde-Alonso has also been supported by the Madrid Government (Comunidad de Madrid-Spain) under the Multiannual Agreement with Universidad Aut\'onoma de Madrid in the line of action encouraging youth research doctors, in the context of the V PRICIT, project SI1/PJI/2019-00514.

\bibliography{BibliographyWeakL1}

\bibliographystyle{siam}

\vskip10pt

\hfill \noindent \textbf{L\'eonard Cadilhac} \\
\null \hfill Laboratoire de Math\'ematiques d'Orsay
\\ \null \hfill Paris-Saclay University
\\ \null \hfill F-91405 Orsay Cedex 
\\ \null \hfill\texttt{leonard.cadilhac@universite-paris-saclay.fr}

\vskip5pt

\hfill \noindent \textbf{Jos\'e M. Conde-Alonso} \\
\null \hfill Departamento de Matem\'aticas
\\ \null \hfill Universidad Aut\'onoma de Madrid 
\\ \null \hfill Cantoblanco. 28049, Madrid. Spain 
\\ \null \hfill\texttt{jose.conde@uam.es}

\vskip5pt

\hfill \noindent \textbf{Javier Parcet} \\
\null \hfill Instituto de Ciencias Matem{\'a}ticas 
\\ \null \hfill Consejo Superior de Investigaciones Cient{\'\i}ficas 
\\ \null \hfill C/ Nicol\'as Cabrera 13-15. 28049, Madrid. Spain 
\\ \null \hfill\texttt{parcet@icmat.es}

\end{document}